\documentclass[11pt]{article}
\usepackage{amsfonts}
\usepackage{color}
\usepackage{graphics}

\usepackage{cite}
\usepackage{latexsym}
\usepackage{amsmath}
\usepackage{amssymb}
\usepackage[dvips]{epsfig}
\usepackage{amscd}

 \usepackage{color}
 \usepackage{indentfirst}
\usepackage{ntheorem}
\theoremstyle{plain}
\theoremseparator{.}
\newtheorem{theorem}{Theorem}[section]
\newtheorem{remark}{Remark}[section]

\newtheorem{lemma}[theorem]{Lemma}

\newcommand\thmref[1]{Theorem~\ref{#1}}
\newcommand\lemref[1]{Lemma~\ref{#1}}

\newcommand{\nb}{{\bold b}}
\newcommand{\nw}{{\bold w}}
\newcommand{\ti}{\tilde}

\newcommand{\lm}{\lambda}

\def\pf{{\it Proof.}  }

\newcommand{\thatsall}{\hfill$\Box$}
\newcommand{\bi}{\bibitem}

\newcommand{\bt}{\begin{theorem}}
\newcommand{\bl}{\begin{lemma}}
\newcommand{\el}{\end{lemma}}
\newcommand{\et}{\end{theorem}}

\renewcommand{\b}{\beta  }

\newcommand{\te}{\theta}

\newcommand{\al}{\alpha}

\newcommand{\ve}{\varepsilon}
\newcommand{\la}{\label}

\newcommand{\ka}{\kappa}

\newcommand{\bn}{\begin{eqnarray}}
\newcommand{\en}{\end{eqnarray}}
\newcommand{\bnn}{\begin{eqnarray*}}
\newcommand{\enn}{\end{eqnarray*}}

\newcommand{\bnnn}{\begin{eqnarray*}}
\newcommand{\ennn}{\end{eqnarray*}}
\newcommand{\ben}{\begin{enumerate}}
\newcommand{\een}{\end{enumerate}}
\newcommand{\ba}{\begin{aligned}}
\newcommand{\ea}{\end{aligned}}
\newcommand{\be}{\begin{equation}}
\newcommand{\ee}{\end{equation}}

\def\norm[#1]#2{\|#2\|_{#1}}

\def\O{\Omega}

\makeatletter
\@addtoreset{equation}{section}
\makeatother

\makeatletter
\@addtoreset{equation}{section}
\makeatother

\title{ On the Global Strong Solutions to   Magnetohydrodynamics  with
Density-Dependent Viscosity and  Degenerate Heat-Conductivity in Unbounded Domains 
}

\author{Yuebo Cao$^1$,  Yi Peng$^2$, Ying Sun$^1$ \thanks{Email addresses:    94099907@qq.com (Y.B.Cao), pengyi16@mails.ucas.ac.cn (Y.Peng), 1913349041@qq.com (Y. Sun)
}
\\[3mm]   1.  School of  Mathematical Sciences,   Xiamen University, \\ Xiamen 361005, P. R. China \\[3mm]   2.   School of Mathematical Sciences, \\ University of Chinese Academy of Sciences,\\ Beijing 100049,  P. R. China}
\date{ }

\begin{document}
\maketitle
\begin{abstract}
 For the equations of a planar magnetohydrodynamic (MHD) compressible flow  with the viscosity  depending  on the specific volume of the gas   and the heat conductivity being proportional to a positive power of the  temperature, we obtain   global existence   of the unique  strong solutions  to the Cauchy problem or the initial-boundary-value one  under natural    conditions   on  the  initial data in one-dimensional unbounded domains.
Our result generalizes the classical one  of the compressible Navier-Stokes system with  constant viscosity and  heat conductivity ([Kazhikhov. Siberian Math. J. (1982)])  to the planar MHD  compressible flow with nonlinear viscosity and degenerate heat-conductivity,  which means   no shock wave, vacuum, or mass or heat concentration will be developed in finite time, although  the  interaction between the magnetodynamic effects and hydrodynamic    is complex and the motion of the flow has large oscillations.
 \end{abstract}

$\mathbf{Keywords.}$ Magnetohydrodynamics,   Strong solutions, Degenerate heat-conductivity, Density-dependent viscosity,  Unbounded Domains

{\bf Math Subject Classification:} 35Q35; 76N10.

\section{Introduction}

Magnetohydrodynamics (MHD), concerning  the dynamics of magnetic fields in electrically conducting fluids, e.g. in plasmas and liquid metals,
covers a wide range of physical objects from liquid metals to cosmic plasmas (\cite{q1,q2,q3,q4,q5,q6,q7}). The central point of MHD theory is that conductive fluids can support magnetic fields. The presence of magnetic fields leads to forces that in turn act on the fluid (typically a plasma), thereby potentially altering the geometry (or topology) and strength of the magnetic fields themselves. 
 We are concerned with the governing equations of a planar magnetohydrodynamic compressible flow written
  in the Lagrange variables
  \be\la{1.1}
v_t=u_{x},
\ee
\be\la{1.2}
u_{t}+(P+\frac12 |\nb|^2)_{x}=\left(\mu\frac{u_{x}}{v}\right)_{x},\ee
 \be\la{1.3}\nw_t-\nb_x=\left(\lm\frac{\nw_{x}}{v}\right)_{x}, \ee \be\la{1.4} (v\nb)_t-\nw_x =\left(\nu\frac{\nb_{x}}{v}\right)_{x}, \ee
\be\la{1.5}\ba&\left(e+\frac{u^2+|\nw|^2+v|\nb|^2}{2}\right)_{t}+ \left(
u\left(P+\frac12|\nb|^2\right)-\nw\cdot\nb\right)_{x}\\&\quad=\left(\kappa\frac{\theta_{x}}{v}
+\mu\frac{uu_x}{v}+\lm\frac{\nw\cdot\nw_x}{v}+\nu\frac{\nb\cdot\nb_x}{v}\right)_x,
 \ea
\ee
where $t>0$ is time, $x\in \Omega\subset \mathbb{R}=(-\infty,+\infty)$ denotes the
Lagrange mass coordinate,  and the unknown functions $v>0, u, \nw\in \mathbb{R}^2, \nb\in   \mathbb{R}^2, e>0, \theta>0$ and $P$ are,  respectively, the specific volume of the gas, longitudinal velocity, transverse velocity, transverse magnetic field, internal energy,  absolute temperature and  pressure. $\mu$ and $\lm$ are the viscosity of the flow,
$\nu$ is the magnetic diffusivity of the magnetic field, and $\ka$ is the heat conductivity.

In this paper, we
concentrate on  a perfect gas for magnetohydrodynamic flow, that is, $P$ and $e$ satisfy \be \la{1.6}   P =R \theta/{v},\quad e=c_v\theta +\mbox{const},
\ee
where  both specific gas constant  $R$ and   heat
capacity at constant volume $c_v $ are   positive constants.
We also assume that $\lm $ and $\nu$ are positive constants, and that $\mu, \ka$ satisfy \be\la{1.7}   \mu=\ti\mu_1+\ti\mu_2 v^{-\al} , \quad \ka=\ti\ka \te^\beta, \ee
with constants $\ti\mu_1>0,\ti \mu_2\ge 0,\ti\ka>0,$  $\al\ge 0, $ and $\beta\ge 0.$

The system \eqref{1.1}-\eqref{1.7} is supplemented
with the initial  conditions
\be\la{1.8}(v,u,\te,\nb,\nw)(x,0)=(v_0,u_0,\te_0,\nb_0,\nw_0)(x),  \quad  x\in\Omega, \ee
and three types of far-field and boundary  ones:

 1) Cauchy problem
\be\la{1.9}\Omega=\mathbb{R}, \quad\lim_{|x|\rightarrow\infty}\left(v,u,\theta,\nb,\nw\right)=(1,0,1,0,0),\quad t>0;\ee

 2) boundary and far-field conditions for $\Omega=(0,\infty)$
\be\ba\la{1.10} &u(0,t)=0,\quad\theta(0,t)=1,\quad\nb(0,t)=\nw(0,t)=0,\,\\&\lim_{x\rightarrow\infty}\left(v,u,\theta,\nb,\nw\right)=(1,0,1,0,0),\quad t>0; \ea\ee

 3) boundary and far-field conditions for $\Omega=(0,\infty)$
\be\ba\la{1.11} &u(0,t)=\theta_x(0,t)=0,\quad\nb(0,t)=\nw(0,t)=0,\\& \lim_{x\rightarrow\infty}\left(v,u,\theta,\nb,\nw\right)=(1,0,1,0,0),\quad t>0. \ea\ee

  There is huge literature on the studies of the global existence   of
 solutions to the  compressible Navier-Stokes system and MHD. Indeed,
for compressible Navier-Stokes system,  Kazhikhov and Shelukhin
\cite{9} first obtained   the global existence of solutions in bounded domains
 for constant coefficients $(\al=\beta=0)$  with large initial data.  From then on,
significant progress has been made on the mathematical aspect of the initial and
initial boundary value problems for $\al=0,\b\ge 0$, see \cite{6,cw1,cw2,fjn1,25,fhl1,24,28,hs1} and the references therein.
 For the Cauchy problem \eqref{1.1}-\eqref{1.9} and the initial boundary value problems \eqref{1.1}-\eqref{1.8}, \eqref{1.10}
and \eqref{1.1}-\eqref{1.8}, \eqref{1.11} (in unbounded domains), Kazhikhov \cite{kazh}  (also cf. \cite{anto}) first for constant coefficients $(\al=\beta=0)$ and very recently, Li-Shu-Xu \cite{py} for $\al=0,\b>0$ obtain
the global existence of strong solutions.

 As for MHD, the existence and uniqueness of
local smooth solutions was first proved in \cite{vh} for bounded domains.   The global existence of strong solutions with large initial data was obtained by
    Kazhikhov  \cite{ka1,az1} for constant coefficients $(\al=\beta=0)$  and   \cite{7,hss1} for  $\al\ge 0 $ and $\b>0$. However, the  methods used in \cite{7,hss1} depend heavily on the boundedness of the domain $\Omega$ and cannot be adapted directly to the case of unbounded domains. Therefore, in this paper, we will prove the global existence of strong solutions to the Cauchy problem \eqref{1.1}-\eqref{1.9}, and the initial boundary value problem \eqref{1.1}-\eqref{1.8} \eqref{1.10} and \eqref{1.1}-\eqref{1.8}\eqref{1.11}. That is, our main result is as follows.

 \begin{theorem}\la{thm1.1} Suppose that \be \la{a1.10}\al\geq0,\quad \beta> 0,\ee
 and that the initial data $ ( v_0,u_0,\te_0,\nb_0,\nw_0)$   satisfies
  \be  \la{a1.11} ( v_0-1,u_0,\te_0-1,\nb_0,\nw_0)\in H^1(\O), \ee  and \be \la{a1.12}
\inf_{x\in \O}v_0(x)>0, \quad \inf_{x\in \O }\theta_0(x)>0, \ee
and are compatible with \eqref{1.10}, \eqref{1.11}. Then there exists a unique global strong solution $(v,u,\theta,\nb,\nw)$  to the initial-boundary-value problem \eqref{1.1}-\eqref{1.9}, or \eqref{1.1}-\eqref{1.8} \eqref{1.10}, or \eqref{1.1}-\eqref{1.8} \eqref{1.11} satisfying for any $T>0$,
 \be
 \begin{cases} \la{a1.13}  v-1, u,\,\theta-1,\nb,\nw \in L^\infty(0,T;H^1(\O)),\\ v_t\in
  L^\infty(0,T;L^2(\O))\cap L^2(0,T;H^1(\O)), \\ u_t,\,\theta_t,\,\nb_t,\,\nw_t,\,u_{xx},\,\te_{xx},\,\nb_{xx},\,\nw_{xx} \in
  L^2((\O)\times(0,T)),\end{cases}\ee
  and for each $(x,t)\in \O \times[0,T]$
  \be C^{-1}\leq v(x,t)\leq C,\quad C^{-1}\leq\te(x,t)\leq C,\ee
  where $C>0$ is a constant  depending on the   data and T.
  \end{theorem}
  A few remarks are in order.
\begin{remark} Our result can be regarded as a natural generalization of   Kazhikhov's theory \cite{kazh}  for the constant viscosity and heat conductivity case to the degenerate and nonlinear one that $\al\geq 0,\beta>0$.\end{remark}
 \begin{remark} Our  result still holds for compressible Navier-Stokes system ($\nb\equiv 0, \nw\equiv 0$) which generalized sightly  those due to \cite{py} where they only consider the case $\al=0,\b>0.$\end{remark}

We now make comments on the analysis of this paper. Compared with  the case of the compressible Navier-Stokes system (\cite{py}),
the main difficulty comes from the nonlinearity of the viscosity and the complex  interaction between the hydrodynamic and magnetodynamic effects in unbounded domains.  The key observations are as follows:
First,   motivated by    \cite{kazh,py}, we obtain an explicit expression of $v$ (see \eqref{2.000}) which is useful for getting the lower and upper bounds of $v.$ Indeed,  on the one hand, observing that there exists a universal constant $\al_1$ such that the integral of $v1_{\{v>\al_1/2\}}$ with respect to $x$ over any interval with length 1 has a lower bound $\al_1/2$ (see \eqref{21q}), we can bound $v$ from below (see \eqref{qp2.17}). On the other hand, observing that the measure of the set $\left\{x\in\O|\te(x,t)>2\right\} $ is uniformly  bounded,  we   multiply \eqref{2.1} by $(\te^\eta-2^\eta)_{+}\te^{\eta-1} $  and use Gr\"onwall's inequality to bound  $v$ from above (see \eqref{2.34}). Next,  multiplying the  equation of the temperature, \eqref{2.1}, by $\te^{-2}(\te^{-1}-2)_{+}^p$ and noticing that the domain $\left\{x\in\O |\theta(x,t)<1/2\right\}$ remains bounded for all $ t\in[0,T]$ (see \eqref{2.18}),  we find that the temperature is indeed bounded from below (see Lemma \ref{lemma2000}), which lays a firm foundation for our further analysis. Finally, to obtain the higher order estimates, we will modify some ideas due to   \cite{py,L-L} to obtain the estimates on the norms of both $u_t$ and $u_{xx}$  (see Lemma \ref{lemma70})  which are crucial for further estimates on the upper bound of the temperature $\te$. The whole procedure will be carried out in the next section.

\section{ Proof of \thmref{thm1.1}}
 We first state   the following  existence and uniqueness of local solutions which can be obtained by using the Banach theorem and the contractivity of the operator defined by the linearization of the problem on a small time interval (c.f. \cite{10,13,tan}).
\begin{lemma} \la{plq1}Let \eqref{a1.10}-\eqref{a1.12} hold. Then there exists some $T>0$ such that  the initial-boundary-value problem \eqref{1.1}-\eqref{1.9}, or \eqref{1.1}-\eqref{1.8} \eqref{1.10}, or \eqref{1.1}-\eqref{1.8} \eqref{1.11} has a unique strong solution $(v,u,\te,\nb,\nw)$ with positive $v(x,t)$ and $\theta(x,t)$ satisfying \eqref{a1.13}.
\end{lemma}

Theorem \ref{thm1.1} will be proved by extending the local solutions globally in time based on the global a priori estimates of solutions (see Lemma \ref{lemma30}--\ref{lemma80}) which will be obtained below.

Without loss of generality, we assume that $\lambda=\nu=\ti\ka=R=c_v=\ti\mu_1=1 , \ti\mu_2=\al.$


\begin{lemma}\la{lemma30} It holds that for any $(x,t)\in\O\times [0,T],$
\be\ba\la{2.12} C^{-1}\le v(x,t)\le C ,\ea\ee
 where (and in what follows)   $C $   denotes some generic positive constant
 depending only on $T, \al, \b,\|(v_0-1,u_0,\theta_0-1 ,\nb_0,\nw_0)\|_{H^1(\O)},
 \inf\limits_{x\in\O} v_0(x),$ and $ \inf\limits_{x\in\O}\theta_0(x).$
\end{lemma}

\pf
First, for any $x\in\O$, denoting $N=[x],$ we write \eqref{1.2} as
\be\ba\la{2.8} u_t=\sigma_x, \ea\ee
where
\bnn\ba\la{2.9}
\sigma&\buildrel\Delta\over=\frac{\mu u_x}{v}-\left(\frac{\theta}{v}+\frac{1}{2}|\nb|^2\right) =(\ln v-v^{-\al})_t-\left(\frac{\theta}{v}+\frac{1}{2}|\nb|^2\right).
\ea\enn
Integrating \eqref{2.8} over $[N,x]\times[0,t]$ leads to
\be\ba
\int_N^{x}u dy-\int_N^{x}u_0 dy=&\ln v- v^{-\al}-\ln v_0+ v_0^{-\al}\\&-\int_0^t\left(\frac{\theta}{v}+\frac{1}{2}|\nb|^2 \right)d\tau
-\int_0^t\sigma(N,\tau)d\tau,\nonumber
\ea\ee
 which implies
\bn\la{2.11}\ba v(x,t) =   B_N(x,t)Y_N(t) \exp\left\{ v^{-\al}\right\} \exp\left\{ \int_0^t\left(\frac{\te}{v}+\frac{1}{2}|\nb|^2\right)d\tau\right\}, \ea\en
where
\be\ba \la{2.0}B_N(x,t)\triangleq v_0\exp\{- v_0^{-\al}\} \exp\left\{\int_N^x u dy-\int_N^x u_0 dy\right\},  \ea\ee
and
\be\ba\la{2.00} Y_N(t)\triangleq\exp\left\{\int_0^t \sigma(N,\tau)d\tau\right\}. \ea\ee
 Denoting
 \be\ba g(x,t)=\int_0^t\left(\frac{\te+\frac{1}{2}v|\nb|^2}{v}\right)d\tau,\nonumber\ea\ee
 we have by \eqref{2.11}
 \be\ba g_t=\frac{\theta+\frac{1}{2}v|\nb|^2}{v}=\frac{\theta+\frac{1}{2}v|\nb|^2}{B_N(x,t)Y_N(t)\exp\{\al v^{-\al}\}\exp\left\{g\right\}},\nonumber\ea\ee
 which gives
 \be\ba \exp\left\{g\right\}=1+\int_0^t\frac{\te+\frac{1}{2}v|\nb|^2}{B_N(x,\tau)Y_N(\tau)\exp\{\al v^{-\al}\}}d\tau .\nonumber \ea\ee
Thus, it follows from  \eqref{2.11}  that
\be \la{2.000}\ba v(x,t)=   B_N (x,t)Y_N (t)\exp\{ v^{-\al}\}\left(1+\int_0^t\frac{\exp\{- v^{-\al}\}(\te+\frac{1}{2}v|\nb|^2)}{B_N(x,\tau)Y_N(\tau)}d\tau \right).\ea\ee

Next, using  \eqref{1.1}-\eqref{1.4}, we rewrite the energy equation \eqref{1.5}   as
\be\la{2.1}\theta_{t}+  \frac{\theta}{v}u_{x}= \left(\frac{\theta^\b\theta_{x}}{v}\right)_{x}+ \frac{ \mu u_{x}^{2}+|\nw_x|^2+|\nb_x|^2}{v}.\ee
Multiplying   \eqref{1.1}, \eqref{1.2}, \eqref{1.3}, \eqref{1.4}, and \eqref{2.1} by $ 1- {v}^{-1} ,
u, \nw, \nb,$ and  $  1- {\theta}^{-1} $  respectively and adding them altogether, we obtain
\bnn\ba &\left(\frac{u^2+|\nw|^2+v|\nb|^2}{2}+\left(\theta-\ln \theta-1\right)+\left(v-\ln v-1\right)\right)_{t}\\&\quad +\frac{\theta^\b \theta_{x}^2}{v\theta^2}+\frac{\mu u_x^2+|\nw_{x}|^2+|\nb_{x}|^2}{v\theta}\\
&=\left(\frac{\theta^{\b}\theta_{x}}{v}+\frac{\mu u u_x}{v}+\frac{\nw\cdot\nw_{x}}{v}+\frac{\nb\cdot\nb_{x}}{v}\right)_{x}+u_x \\&\quad-\left(u\left(\frac{\theta}{v}+\frac{1}{2}|\nb|^2\right)
-\nw\cdot\nb\right)_{x}-\left(\frac{\theta^{\b}\theta_{x}}{v\theta}\right)_x,
\ea\enn
which together with \eqref{1.9} or \eqref{1.10} or \eqref{1.11} yields
\be\ba\la{2.2}&\sup_{0\le t\le T}\int_\O\left(\frac{u^2+|\nw|^2+v|\nb|^2}{2}+(v-\ln v-1) +(\theta-\ln \theta-1)\right)dx
\\&\quad+\int_0^T W(t)dt\le e_0,      \ea\ee
where \be\la{2.3}W(t)\buildrel \Delta \over =\int_\O\left(\frac{\theta^\b\theta_{x}^2}{v\theta^2}+\frac{\mu u_x^2+|\nw_x|^2+|\nb_x|^2}{v\theta}\right)dx, \nonumber\ee
and
\be e_{0}\triangleq  2\int_\O\left(\frac{ u_0^2+|\nw_{0}|^2+v_0|\nb_{0}|^2}{2}+(v_0-\ln v_0-1)+(\theta_0-\ln \theta_0-1)\right)dx\nonumber . \ee

Next,  we have by \eqref{2.2}
$$ \int_{N}^{N+1} (v-\ln v-1+\theta-\ln \theta-1)dx \le e_{0},$$
which together with Jensen's inequality yields that for any $t\in [0,T],$
\be\la{2.5}\al_1\le\int_N^{N+1}v(x,t)dx\le\al_2,\quad\al_1\le\int_N^{N+1}\theta(x,t)dx\le\al_2,
\ee
where $0<\al_1<\al_2$ are two roots of
\be z-\ln z-1=e_0.\nonumber\ee
Moreover, it follows from \eqref{2.2}  that
\be\ba\left|\int_N^x\left(u(y,t)-u_0(y)\right)dy\right|\le\left(\int_N^{N+1} u^2 dy\right)^{\frac{1}{2}}+\left(\int_N^{N+1}u_0^2dy\right)^{\frac{1}{2}}\le C,\nonumber\ea\ee
which together with \eqref{2.0} implies
\be\la{2.13}    C^{-1}\le B_N(x,t)\le C,  \ee
where and in what follows, $C$ is a constant independent of $N$.

Next, letting
\bnn 1_{\left\{v >\frac{\al_1}{2}\right\}}=\begin{cases} 1,\, &\mbox{if} \,\, v >\frac{\al_1}{2}, \\  0, &\mbox{if} \,\, v \le\frac{\al_1}{2} , \end{cases} \enn
multiplying \eqref{2.000} by $\frac{1}{Y_N(t)}1_{\left\{v >\frac{\al_1}{2}\right\}} $ and  integrating the resultant equality over $[N,N+1]$ leads to
 \be\ba\la{2.14} &\frac{1}{Y_N(t)}\int_N^{N+1}v(x,t)1_{\left\{v >\frac{\al_1}{2}\right\}}dx\\&\le C\int_N^{N+1}B_N(x,t)\exp\left\{
 v^{-\al}\right\}1_{\left\{v >\frac{\al_1}{2}\right\}}\left(1+\int_0^t\frac{\te+\frac{1}{2}v|\nb|^2}{B_N(x,\tau)Y_N(\tau)}d\tau\right) dx \\&
\le C +C\int_0^t\frac{1}{Y_N(\tau)}\int_N^{N+1}\left(\theta+\frac{1}{2}v|\nb|^2\right)dxd\tau
 \\&\le C+C\int_0^t\frac{1}{Y_N(\tau)}d\tau,\ea\ee
 where we have used the \eqref{2.2} and \eqref{2.13}. Since  \eqref{2.5} implies
 \be \la{21q}\ba     \int_N^{N+1} v(x,t) 1_{\left\{v >\frac{\al_1}{2}\right\}}dx\ge \frac{\al_1}{2} , \ea\ee
 it follows from  \eqref{2.14}   and  Gr\"onwall's inequality   that for any $t\in [0,T],$
\be\la{2.15} \frac{1}{Y_N(t)} \le C.\ee

Next, integrating \eqref{2.000}  over $[N,N+1]$, after using \eqref{2.13} we obtain that
\be\ba \int_N^{N+1}vdx &\geq Y_N(t)\int_N^{N+1}B_N(x,t)\exp\left\{ v^{-\al}\right\}dx
\\&\geq Y_N(t)\int_N^{N+1}B_N(x,t)dx\\&\geq C^{-1}Y_N(t),
  \nonumber\ea\ee
which together with \eqref{2.5}  yields  that for any $t\in [0,T],$
\bnn Y_N(t) \le C .  \enn
Combining this, \eqref{2.000}, \eqref{2.15},   and \eqref{2.13}  gives that for $(x,t)\in [N,N+1]\times [0,T]$
\be\ba C^{-1}\le v(x,t)\le C+C\int_0^T\left(\sup_{x\in[N,N+1]}\te+\sup_{x\in[N,N+1]}v
\sup_{x\in[N,N+1]}| \nb|^2\right)dt, \nonumber \ea\ee
which in particular implies that for any $(x,t)\in\O\times [0,T]$
\be\la{qp2.17}\ba C^{-1}\le v(x,t)\le C+C\int_0^T\left(\sup_{x\in\O}\te+\sup_{x\in\O}v\sup_{x\in\O}|\nb|^2\right)dt.  \ea\ee

Next, denoting   \be\ba\left(\theta >2\right)(t)=\left\{x\in\O|\theta(x,t)>2\right\}, \nonumber\ea\ee
and \be\left(\theta<1/2\right)(t)=\left\{x\in\O |\theta(x,t)<1/2\right\}.\nonumber\ee
we get by  \eqref{2.2}
\be\ba e_0&\geq\int_{\left(\te<1/2\right)(t)}\left(\te-\ln \te-1\right) dx+\int_{\left(\te>2\right)(t)}\left(\te-\ln \te-1\right)dx\\&
\geq\left(\ln 2-1/2\right)\left|\left(\te< 1/2\right)(t)\right|+\left(1-\ln 2\right)\left|\left(\te>2\right)(t)\right|\\&
\geq\left(\ln 2-1/2\right)\left(|(\te<1/2)(t)|+|\left(\te>2\right)(t)|\right), \nonumber\ea\ee
which shows that for any $t\in[0,T]$
\be\ba\la{2.18}|(\te< 1/2)(t)|+|(\te>2)(t)|\le\frac{2e_0}{2\ln2-1}.  \ea\ee
For    $\eta\triangleq \max\left\{1, {2-\b} \right\} /4\in\left(0,\frac{1}{2}\right)$, integrating \eqref{2.1} multiplied by
 $\left(\te^{\eta}-2^{\eta}\right)_{+}\te^{\eta-1}$ over $\O\times(0,T)$, we get
 \be\ba\la{2.24}&(1-2\eta)\int_0^T\int_{(\te>2)(t)} \frac{\te^\b\te_x^2}{v\te^{2-2\eta}}dxdt\\&\quad+\int_0^T\int_{\O}\frac{\mu u_x^2+|\nw_x|^2+|\nb_x|^2}{v}
 \left(\te^\eta-2^\eta\right)_{+}\te^{\eta-1} dxdt\\&
 =\frac{1}{2\eta}\int_{\O}\left(\left(\te^\eta-2^\eta\right)_{+}^2-\left(\te_0^\eta-2^\eta\right)_{+}^2\right)dx+2^\eta(1-\eta)\int_0^T\int_{(\te>2)(t)}\frac{\te^\b\te_x^2}{v\te^{2-\eta}}dxdt\\&\quad
 +\int_0^T\int_{\O}\frac{\te u_x}{v}\left(\te^\eta-2^\eta\right)_{+}\te^{\eta-1}dxdt\\&
 \le C+\frac{1-2\eta}{2}\int_0^T \int_{(\te>2)(t)} \frac{\te^\b\te_x^2}{v\te^{2-2\eta}}dxdt+\frac{1}{2}\int_0^T\int_{\O}\frac{\mu u_x^2}{v}\left(\te^\eta-2^\eta\right)_{+}\te^{\eta-1}dxdt\\&\quad+C\int_0^T
 \int_{\O}\frac{\te^2}{\mu v}\left(\te^\eta-2^\eta\right)_{+}\te^{\eta-1}dxdt,   \ea\ee
where in the last inequality we have used \eqref{2.2} and the following inequality
\be\ba\la{2.25}\sup_{0\le t\le T}\int_{(\te>2)(t)}\te dx\le C\sup_{0\le t\le T}\int_{\O} (\te-\ln \te-1)dx\le C. \ea\ee
Using \eqref{2.25}, \eqref{2.18}, \eqref{qp2.17}, and Young's inequality, we get
\be\ba\la{2.26}&\int_0^T\int_{\O}\frac{1}{\mu v}\left(\te^\eta-2^\eta\right)_{+}\te^{\eta+1} dxdt\\&
\le C\int_0^T\int_{(\te>2)(t)}\left(\te^\eta-2^\eta\right)\te^{\eta+1}dxdt\\&
\le C\int_0^T\sup_{x\in\O}\te^{2\eta}\int_{(\te>2)(t)}\te dxdt+ C\int_0^T\sup_{x\in\O}\te^{2\eta}\int_{(\te>2)(t)}\te^{1-\eta} dxdt\\&
\le C\int_0^T\sup_{x\in{\O}}\te^{2\eta}dt\\&
\le \varepsilon\int_0^T\sup_{x\in{\O}}\te dt+C(\varepsilon). \ea\ee

We then deduce from Cauchy's inequality that
\be\ba\la{2.27} \int_0^T\sup_{x\in{\O}}\te dt  &
\le C\int_0^T\sup_{x\in\O}\left(\int_{\infty}^x\partial_{y}(\te-2)_{+}(y,t)dy \right)dt+C \\&
\le C\int_0^T\int_{(\te>2)(t)}\frac{\te^\b \te_x^2}{v\te^{2-2\eta}}dxdt+C\int_0^T\int_{(\te>2)(t)}\frac{v\te^{2-2\eta}}{\te^\b} dxdt+C \\&\le C\int_0^T\int_{(\te>2)(t)}\frac{\te^\b\te_x^2}{v\te^{2-2\eta}} dxdt +C\\&\quad +C\int_0^T\sup_{x\in\O}\te^{\max\{2-2\eta-\b,0\}} \int_{ (\te>2)(t) } v dxdt   \\&
\le C\int_0^T \int_{(\te>2)(t)}\frac{\te^\b \te_x^2}{v\te^{2-2\eta}}dxdt+\ve\int_0^T\sup_{x\in\O}\te dt+C(\ve), \ea\ee
where in the last inequality we have used \eqref{2.18} and the follow inequality
\be\ba \sup_{0\le t\le T}\int_{(\te>2)(t)}v dx&\le \sup_{0\le t\le T}\int_{(\te>2)(t)\cap (v\le 2)(t)}v dx+\sup_{0\le t\le T}\int_{(v> 2)(t)}v dx \\&\le C+ C\sup_{0\le t\le T}\int_{\O}\left(v-\ln v-1\right)dx\le C.  \nonumber \ea\ee
Putting \eqref{2.26} and \eqref{2.27} into \eqref{2.24} gives
\be\ba \la{400}\int_0^T\sup_{x\in\O}\te dt+\int_0^T\int_{(\te>2)(t)}\frac{\te^\b\te_x^2}{v\te^{2-2\eta}}dxdt\le C, \ea\ee  which together with   \eqref{2.2} yields
\be\ba\la{2.300} \int_0^T\int_{\O}\frac{\te^\b\te_x^2}{v\te^{2-2\eta}}dxdt \le C +\int_0^T\int_{\left(\te\le 2\right)(t)}\frac{\te^\b\te_x^2}{v\te^{2}}dxdt\le C. \ea\ee

Finally, using \eqref{2.2} and \eqref{400}, we have
\be\ba\la{2.22}\int_0^T\sup_{x\in{\O}}|\nb|^2dt& \le C\int_0^T\int_{\O}|\nb\cdot\nb_x|dxdt \\&
\le C\int_0^T\int_{\O}\frac{|\nb_x|^2}{v\te}dxdt+C\int_0^T \int_{\O}\te v |\nb|^2dxdt\\&
\le C+C\int_0^T\sup_{x\in\O}\te dt\\&
\le C,\ea\ee
which together with \eqref{qp2.17}, \eqref{400}, \eqref{2.22} and Gr\"onwall's inequality yields
\be\la{2.34}\sup_{(x,t)\in\O\times [0,T]} v\le C.   \ee
Combining this with \eqref{2.12} finishes the proof of \lemref{lemma30}.  \thatsall
\begin{lemma}\la{lemma2000}
There exists a positive constant $C$ such that for all $(x,t)\in\O\times[0,T],$
\be\la{2.q27}\ba \theta(x,t)\geq C^{-1}. \ea\ee
\end{lemma}

\pf For any $p>2$,  multiplying \eqref{2.1} by $\te^{-2}\left(\te^{-1}-2\right)_{+}^{p}$ with $\left(\te^{-1}-2\right )_{+}\buildrel \Delta \over =\max\left\{\te^{-1}-2,0\right\}$  and integrating the resultant equality over $\O$, we obtain
\bnn\ba &\frac{1}{p+1}\frac{d}{dt}\int_{\O}\left(\frac{1}{\te}-2\right)_{+}^{p+1}dx+\int_{\O}\frac{\mu u_x^2+|\nw_x|^2+|\nb_x|^2}{v\te^2}\left(\frac{1}{\te}-2\right)_{+}^pdx \\&\le\int_{\O}\frac{u_x}{v\te}\left(\frac{1}{\te}-2\right)_{+}^p dx\\&
\le\frac{1}{2}\int_{\O}\frac{\mu u_x^2}{v\te^2}\left(\frac{1}{\te}-2\right)_{+}^pdx+\frac{1}{2}\int_{\O}\frac{1}{\mu v}
\left(\frac{1}{\te}-2\right)_{+}^p dx \\&
\le\frac{1}{2}\int_{\O}\frac{\mu u_x^2}{v\te^2}\left(\frac{1}{\te}-2\right)_{+}^pdx+C\int_{\O}\left(\frac{1}{\te}-2\right)_{+}^pdx,   \ea\enn
where in the last inequality we have used \eqref{2.12}. Thus, combining this with  \eqref{2.18} leads to
\be\ba&\left\|\left(\frac{1}{\te}-2\right)_{+}\right\|_{L^{p+1}(\O)}^{p} \frac{d}{dt}\left\|\left(\frac{1}{\te}-2\right)_{+} \right\|_{L^{p+1}(\O)}\\& \le C \int_{\O}\left(\frac{1}{\te}-2\right)_{+}^p dx \\&
\le C\left(\int_{\O}\left(\frac{1}{\te}-2\right)_{+}^{p+1}dx\right)^{\frac{p}{p+1}}, \nonumber  \ea\ee
with $C$ independent of $p$. This in particular implies that there exists some positive constant $C$ independent of $p$
such that
\bnn\sup_{0\le t\le T}\left\|\left(\frac{1}{\te}-2\right)_{+}\right\|_{L^{p+1}(\O)}\le C.  \enn
 Letting $p\rightarrow +\infty$  and using \eqref{2.18}   shows
\be\ba\sup_{0\le t\le T}\left\|\left(\frac{1}{\te}-2\right)_{+}\right\|_{L^{\infty}(\O)} \le C, \nonumber\ea\ee
which gives \eqref{2.q27} and finishes the proof of Lemma \ref{lemma2000}. \thatsall
\begin{lemma}\la{lemma40}
There exists a positive constant $C$ such that
\be\ba\la{3.1} \int_0^T\int_{\O}\left(u_x^2+\te^{-1}\te_x^2\right)dxdt\le C.\ea\ee
\end{lemma}

\pf On the one hand, integrating the momentum equality \eqref{1.2} multiplied by $u$ with respect to $x$ over $\O$, we obtain
\be\ba\la{3.00} &\frac{1}{2}\frac{d}{dt}\int_{\O}u^2 dx+\int_{\O}\frac{\mu u_x^2}{v}dx\\&
=\int_{\O}\frac{\te-1}{v}u_x dx-\int_{\O}\frac{(v-1)u_x}{v}dx+\frac{1}{2}\int_{\O}|\nb|^2u_x dx\\&
\le C \int_{\O}(\te-1)^2dx+C\int_{\O}(v-1)^2dx+\frac{1}{2}\int_{\O}\frac{\mu u_x^2}{v}dx+C\int_{\O}|\nb|^4dx\\&
\le C+C\int_{(\te>2)(t)}\te^2dx+C\sup_{x\in\O}|\nb|^2+ \frac{1}{2}\int_{\O}\frac{\mu u_x^2}{v}dx\\&
\le C+C\sup_{x\in\O}\te+C\sup_{x\in \O}|\nb|^2+\frac{1}{2}\int_{\O}\frac{\mu u_x^2}{v}dx,   \ea\ee
where we have used \eqref{2.2}, \eqref{2.25}, and \eqref{2.12}. Integrating \eqref{3.00} in $t$ over $[0,T]$, we get after using \eqref{400} and \eqref{2.22}
\be\ba\la{3.2}\int_0^T\int_{\O}u_x^2dxdt\le C.\ea\ee

On the other hand, we deduce from \eqref{2.300} and \eqref{2.q27} that
\bnn\ba\la{3.4}\int_0^T\int_{\O}\te^{-1}\te_x^2 dxdt&\le C\int_0^T\int_{\O}\frac{\te^\b\te_x^2}{v\te^{2-2\eta}}\cdot\te^{1-2\eta-\b}dxdt\\&
\le C,  \ea\enn  which together with
  \eqref{3.2} finishes the proof of the Lemma \ref{lemma40}.\thatsall
\begin{lemma}\la{lemma50}
There exists a positive constant $C$ such that
\be\ba\la{4.1} \sup_{0\le t\le T}\int_{\O}v_x^2dx\le C.\ea\ee
\end{lemma}

\pf We rewrite the momentum equation \eqref{1.2} as
\bnn\la{4.2}\ba \left(u-\frac{\mu v_x}{v}\right)_t=-\left(\frac{\te}{v}+\frac{1}{2}|\nb|^2\right)_x.   \ea\enn
Multiplying the above equation by $u-\frac{\mu v_x}{v}$ and integrating the resultant equality  yield that for any $t\in(0,T)$
\be\ba\la{4.3}&\frac{1}{2}\int_{\O}\left(u-\frac{\mu v_x}{v}\right)^2 dx -\frac{1}{2}\int_{\O}\left(u-\frac{\mu v_x}{v}\right)(x,0)dx\\&= \int_0^t\int_{\O}\left( \frac{\te v_x}{v^2}-\frac{\te_x}{v}-\nb\cdot\nb_x\right)\left(u-\frac{\mu v_x}{v}\right)dxd\tau \\
&=  -\int_0^t\int_{\O}\frac{ \mu \te v_x^2}{v^3}dxd\tau+\int_0^t\int_{\O}\frac{\te u v_x}{v^2}dxd\tau \\&\quad-\int_0^t\int_{\O}\frac{\te_x}{v} \left(u-\frac{\mu v_x}{v}\right)dxd\tau-\int_0^t\int_{\O}\nb\cdot\nb_x  \left(u-\frac{ \mu v_x}{v}\right)dxd\tau \\
&=  -\int_0^t\int_{\O}\frac{ \mu \te v_x^2}{v^3}dxd\tau+\sum_{i=1}^3I_i.\ea\ee
Each $I_i (i=1,2,3)$ can be estimated as follows:

First, Cauchy's inequality gives
\be\ba\la{4.4}|I_1|&\leq\frac{1}{2}\int_0^t\int_{\O}\frac{\mu\te v_x^2}{v^3}dxd\tau+\frac{1}{2}\int_0^t\int_{\O}\frac{u^2\te}{\mu v}dxd\tau \\&\leq\frac{1}{2}\int_0^t\int_{\O}\frac{\mu\te v_x^2}{v^3}dxd\tau+C\int_0^T\sup_{x\in{\O}}\te d\tau \\&\leq C+\frac{1}{2}\int_0^t\int_{\O}\frac{\mu\te v_x^2}{v^3}dxd\tau,\ea\ee
where  we have used \eqref{2.2}, \eqref{400}, and \eqref{2.12}.

Next, using \eqref{3.1}  and Cauchy's inequality, we have
\be\ba\la{4.5}|I_2| &\leq \frac{1}{2}\int_0^t\int_{\O}\te^{-1}\te_x^2dxd\tau+\frac{1}{2}\int_0^t\int_{\O}\frac{\te}{v^2}\left(u-\frac{\mu v_x}{v}\right)^2dxd\tau \\
&\leq C+C\int_0^t\sup_{x\in{\O}}\te\int_{\O}\left(u-\frac{\mu v_x}{v}\right)^2dxd\tau.\ea\ee

Finally, multiplying \eqref{1.3} by $\nw$, \eqref{1.4} by $\nb$, adding them and integrating the resultant equality in $x$ over $\O$, one has
\bnn\ba&\frac{1}{2}\frac{d}{dt}\int_{\O}\left(v|\nb|^2  + |\nw|^2\right)dx+\int_{\O}\frac{|\nb_x|^2}{v}dx+\int_{\O}\frac{|\nw_x|^2}{v}dx\\&
=-\frac{1}{2}\int_{\O}u_x |\nb|^2dx\\&
\le C\int_{\O}u_x^2dx+C\sup_{x\in\O}|\nb|^2\int_{\O}|\nb|^2dx.
\ea\enn
Integrating this inequality in $t$ over $(0,T)$ and using \eqref{2.2}, \eqref{2.22}, \eqref{2.12} and \eqref{3.1}, we get
\be\ba\la{4.6}\int_0^T\int_{\O}|\nb_x|^2dxdt+\int_0^T\int_{\O}|\nw_x|^2dxdt\le C.\ea\ee
Combining \eqref{4.6} with Cauchy's inequality leads to
\be\ba \la{4.7}|I_3|&\le \frac{1}{2}\int_0^t\int_{\O}|\nb_x|^2dxd\tau+\frac{1}{2}\int_0^t\int_{\O}|\nb|^2\left(u-\frac{\mu v_x}{v}\right)^2dxd\tau \\&
\le C+\frac{1}{2}\int_0^t\sup_{x\in\O}|\nb|^2\int_{\O}\left(u-\frac{\mu v_x}{v}\right)^2dxd\tau.
\ea\ee
Putting \eqref{4.4}, \eqref{4.5}, \eqref{4.7} into \eqref{4.3}, we obtain after using Gr\"onwall's inequality, \eqref{400}, and \eqref{2.22} that
\bnn\ba \sup_{0\le t\le T}\int_{\O}\left(u-\frac{\mu v_x}{v}\right)^2dx+\int_0^T\int_{\O}\frac{\mu\te v_x^2}{v^3}dxdt\leq C,\ea\enn
which together with \eqref{2.2}  gives \eqref{4.1} and finishes the proof of Lemma \ref{lemma50}. \thatsall

\begin{lemma}\la{lemma60}
There is a positive constant C such that
\be \la{6.1}\ba
&\sup_{0\le t\le T} \int_{\O}\left( |\nb_x|^2+|\nw_x|^2\right)dx \\&  +\int_0^T \int_{\O}\left(|\nb_t|^2 +|\nb_{xx}|^2 +|\nw_t|^2+|\nw_{xx}|^2\right)dx dt\leq C.\ea\ee
\end{lemma}

\pf First, multiplying \eqref{1.3} by $\nw_{xx}$ and integrating the resultant equality over $\O\times(0,T)$, we obtain after using \eqref{4.1}, \eqref{4.6} and Cauchy's inequality that
\be\ba\la{6.2}&\frac{1}{2} \int_{\O}|\nw_x|^2dx+\int_0^T \int_{\O}\frac{|\nw_{xx}|^2}{v}dxdt   \\& \leq C+\frac{1}{2}\int_0^T \int_{\O}\frac{|\nw_{xx}|^2}{v}dxdt +C\int_0^T \int_{\O}\left(|\nb_x|^2+|\nw_x|^2v_x^2\right)dxdt \\&\leq C+\frac{1}{2}\int_0^T \int_{\O}\frac{|\nw_{xx}|^2}{v}dxdt +C\int_0^T\sup_{x\in\O}|\nw_x|^2dt.\ea\ee
Direct computation shows  after using \eqref{4.6}
\be\ba\la{6.3}\int_0^T\sup_{x\in\O}|\nw_x|^2dt&\leq C(\ve)\int_0^T \int_{\O}|\nw_x|^2dxdt+\varepsilon\int_0^T \int_{\O}\frac{|\nw_{xx}|^2}{v}dxdt \\&\leq C(\ve)+\varepsilon\int_0^T \int_{\O}\frac{|\nw_{xx}|^2}{v}dxdt,\ea\ee
which combined with \eqref{6.2}  leads to
\bn\ba \la{6.4}\sup_{0\le t\le T} \int_{\O}|\nw_x|^2dx+\int_0^T \int_{\O}|\nw_{xx}|^2dxdt\leq C.\ea\en

Then, we rewrite  \eqref{1.3} as
\bnn\la{6.5}\ba \nw_t=\frac{\nw_{xx}}{v}-\frac{\nw_xv_x}{v^2}+\nb_x, \ea\enn
 which together with \eqref{2.12}, \eqref{4.6}, \eqref{6.3}, \eqref{6.4}, and \eqref{4.1} gives
\be\ba\la{6.6}\int_0^T \int_{\O}|\nw_t|^2dxdt &\leq
C\int_0^T \int_{\O} \left(|\nb_x|^2+|\nw_{xx}|^2 + v_x^2|\nw_x|^2\right)dxdt\\&\leq C\int_0^T\sup_{x\in\O}|\nw_x|^2dt+C\\&\leq C .\ea\ee

Next, multiplying \eqref{1.4} by $\frac{\nb_{xx}}{v}$ and integrating the result over $\O\times(0,T)$, we deduce from \eqref{4.1},  \eqref{3.1}, \eqref{6.4}, \eqref{4.6}, \eqref{2.2}  and Cauchy's inequality that
\bnn\ba &\frac{1}{2}\int_{\O}|\nb_x|^2dx+\int_0^T\int_{\O}\frac{|\nb_{xx}|^2}{v^2}dxdt\\ &\leq C+\frac{1}{2}\int_0^T\int_{\O}\frac{|\nb_{xx}|^2}{v^2}dxdt +C\int_0^T\int_{\O}\left(|\nb_x|^2v_x^2+u_x^2|\nb|^2+|\nw_x|^2\right)dxdt \\&\leq
C+\frac{1}{2}\int_0^T\int_{\O}\frac{|\nb_{xx}|^2}{v^2}dxdt+C\int_0^T\sup_{x\in{\O}}|\nb_x|^2dt +\sup_{(x,t)\in{\O}\times[0,T]}|\nb|^2 \\&\leq C+\frac{3}{4}\int_0^T\int_{\O}\frac{|\nb_{xx}|^2}{v^2}dxdt+C \int_0^T\int_{\O}|\nb_x|^2dxdt \\&\quad  +C\sup_{0\le t\le T}\int_{\O}|\nb|^2dx +\frac{1}{4}\sup_{0\le t\le T}\int_{\O}|\nb_x|^2dx \\&\leq C +\frac{3}{4}\int_0^T\int_{\O}\frac{|\nb_{xx}|^2}{v^2}dxdt+\frac{1}{4}\sup_{0\le t\le T}\int_{\O}|\nb_x|^2dx,\ea\enn which implies
\bn\ba\la{6.7}\sup_{0\le t\le T}\int_{\O}|\nb_x|^2dx+\int_0^T\int_{\O}|\nb_{xx}|^2dxdt\leq C.\ea\en
Hence, \be\la{6.80} \sup_{(x,t)\in{\O}\times [0,T]}|\nb|^2\le C\sup_{0\le t\le T}\int_{\O}|\nb|^2dx+C \sup_{0\le t\le T}\int_{\O}|\nb_x|^2dx\le C.\ee

Finally, we rewrite \eqref{1.4} as
\bnn\la{6.8}\ba\nb_t=\frac{\nw_x}{v}+\frac{\nb_{xx}}{v^2}-\frac{\nb_x v_x}{v^3}-\frac{\nb u_x}{v},\ea\enn
which together with  \eqref{4.1}, \eqref{6.7}, \eqref{6.80} and  \eqref{3.1}  gives
\be\ba\la{6.9}\int_0^T\int_{\O}|\nb_t|^2dxdt&\leq C\int_0^T\int_{\O}\left(|\nb_{xx}|^2+|\nb_x|^2v_x^2+|\nw_x|^2+|\nb|^2u_x^2\right)dxdt\\&\leq C+C\int_0^T\left(\sup_{x\in\O}|\nb_x|^2+ \int_{\O}u_x^2dx\right)dt\\&\leq C+C\int_0^T\int_{\O}\left(|\nb_x|^2+|\nb_{xx}|^2\right)dxdt\\&\leq C.\ea\ee
 Combining \eqref{6.9}, \eqref{6.4}, \eqref{6.6} and  \eqref{6.7}   gives \eqref{6.1} and we finish the proof of Lemma \ref{lemma60}. \thatsall
\begin{lemma}\la{lemma70}
There is a positive constant C such that
\be \la{7.1}\ba
&\sup_{0\le t\le T} \int_{\O} u_x^2 dx   +\int_0^T\int_{\O}\left( u_t^2 +u_{xx}^2 \right)dx dt\leq C.\ea\ee
\end{lemma}

\pf
First, multiplying \eqref{1.2} by $u_{xx}$ and integrating the result over $\O\times(0,T)$,  we have
 \be\ba\la{7.2} &\frac{1}{2}\int_{\O}u_x^2dx+\int_0^T\int_{\O}\frac{\mu u_{xx}^2}{v}dxdt  \\&\leq C+\frac{1}{2}\int_0^T\int_{\O}\frac{\mu u_{xx}^2}{v}dxdt+C\int_0^T\int_{\O}\left(\te_x^2+\te^2v_x^2+|\nb|^2|\nb_x|^2+u_x^2v_x^2 \right)dxdt \\&\leq C+\frac{1}{2}\int_0^T\int_{\O}\frac{\mu u_{xx}^2}{v}dxdt+C\int_0^T\int_{\O}\te_x^2dxdt +C\int_0^T\sup_{x\in\O}\te^2\int_{\O}v_x^2dxdt \\&\quad
+C\sup_{(x,t)\in\O\times[0,T]}|\nb|^2\int_0^T\int_{\O}|\nb_x|^2dxdt
+C\int_0^T\sup_{x\in\O}u_x^2\int_{\O}v^2_xdxdt  \\&\leq C+ \frac{3}{4} \int_0^T\int_{\O}\frac{\mu u_{xx}^2}{v}dxdt
+C\int_0^T\int_{\O}\frac{\te^\b\te_x^2}{v}dxdt+C\int_0^T\sup_{x\in\O}(\te-2)_{+}^2dt\\&
\le C+ \frac{3}{4} \int_0^T\int_{\O}\frac{\mu u_{xx}^2}{v}dxdt
+C_1\int_0^T\int_{\O}\frac{\te^\b\te_x^2}{v}dxdt,
 \ea\ee
where in the third inequality we have used \eqref{4.1}, \eqref{6.1}, \eqref{6.80},  \eqref{3.2}  and the following  inequality
\bnn\ba\int_0^T\sup_{x\in\O}u_x^2dt&\le\int_0^T\int_{\O}|\left(u_x^2\right)_x| dxdt\\&
\le 2\int_0^T\left(\int_{\O}u_{x}^2dx\right)^{\frac{1}{2}}\left(\int_{\O}u_{xx}^2dx\right)^{\frac{1}{2}}dt\\&
\le\frac{1}{4}\int_0^T\int_{\O}\frac{u_{xx}^2}{v}dxdt+C\int_0^T\int_{\O} u_x^2dxdt,\ea\enn
and in the last  inequality we have used
\be\ba\sup_{x\in\O}(\te-2)_{+}^2&=\sup_{x\in\O}\left(\int_x^{\infty}\partial_y(\te-2)_{+}(y,t)dy\right)^2\\&
\le\left(\int_{(\te>2)(t)}|\te_y|dy\right)^2\\&
\le C\int_{\O}\te_x^2dx,   \nonumber\ea\ee due to \eqref{2.18}.

Then, motivated by \cite{L-L}, we integrate \eqref{2.1} multiplied by $(\te-2)_{+}\buildrel\Delta\over=\max\limits_{x\in\O}\left\{\te-2,0\right\}$ over $\O\times(0,T)$ to get

\be\ba\la{7.3}&\frac{1}{2}\int_{\O}(\te-2)_{+}^2dx+\int_0^T\int_{(\te>2)(t)}\frac{\te^\b\te_x^2}{v}dxdt\\&= \frac{1}{2}\int_{\O}(\te_{0}-2)_{+}^2dx+\int_0^T\int_{\O}\frac{\mu u_x^2 + |\nw_x|^2 + |\nb_x^2|}{v}(\te-2)_{+}dxdt\\&\quad-
\int_0^T\int_{\O}\frac{\te(\te-2)_{+}}{v}u_xdxdt\\&
\le C+C\int_0^T\sup_{x\in\O}\te\int_{\O}\left(|\nb_x|^2+|\nw_x|^2\right)dxdt\\&
\quad + C\int_0^T\sup_{x\in\O}\te\left(\int_{\O}(\te-2)_{+}^2dx+\int_{\O}u_x^2dx\right)dt\\&
\le C+C\int_0^T\sup_{x\in\O}\te\left(\int_{\O}(\te-2)_{+}^2dx+\int_{\O}u_x^2dx\right)dt,
\ea\ee
 where we have used \eqref{6.1} and \eqref{400}.
We deduce from \eqref{2.2} that
\be\ba\int_0^T\int_{\O}\te^\b\te_x^2dxdt&=\int_0^T\int_{(\te>2)(t)}\te^\b\te_x^2dxdt+\int_0^T\int_{(\te\le 2)(t)}\te^\b\te_x^2dxdt\\&
\le C\int_0^T\int_{(\te>2)(t)}\frac{\te^\b\te_x^2}{v} dxdt+ C\int_0^T\int_{(\te\le 2)(t)}\frac{\te^\b\te_x^2}{v\te^2} dxdt\\&
\le C\int_0^T\int_{(\te>2)(t)}\frac{\te^\b\te_x^2}{v} dxdt+ C,\nonumber\ea\ee
which together with \eqref{7.3} yields
 \be\ba\la{7.4}&\frac{1}{2}\int_{\O}(\te-2)_{+}^2dx+C_2\int_0^T\int_{\O}\frac{\te^\b\te_x^2}{v}dxdt\\&
\le  C+C\int_0^T\sup_{x\in\O}\te\left(\int_{\O}(\te-2)_{+}^2dx+\int_{\O}u_x^2dx\right)dt.
\ea\ee
Adding \eqref{7.4} multiplied by $2C_2^{-1}C_1$ to \eqref{7.2}, we use Gr\"onwall's inequality and \eqref{400} to get
\be\ba\la{7.5} \sup_{0\le t\le T}\int_{\O}\left(u_x^2+(\te-2)_{+}^2\right)dx+\int_0^T\int_{\O}\left(u_{xx}^2+\te^\b\te_x^2\right)dxdt\le C.  \ea\ee

Finally, we rewrite \eqref{1.2} as
\bnn\ba\la{eq20} u_t=\frac{\mu u_{xx}}{v}+\left(\frac{\mu}{v}\right)^{'}_{v}u_xv_x-\frac{\te_x}{v}+\frac{\te v_x}{v^2}-\nb\cdot\nb_x ,\ea\enn
which together with  \eqref{7.5}, \eqref{4.1} and \eqref{6.1} leads to
\bnn\ba
\int_0^T\int_{\O}u_t^2dxdt&\leq C\int_0^T\int_{\O}\left(u_{xx}^2+u_x^2v_x^2+\te_x^2+\te^2v_x^2+|\nb|^2|\nb_x|^2\right)dxdt\\&\leq C.
\ea\enn
Combining this with  \eqref{7.5}  immediately gives \eqref{7.1} and  completes the proof of \lemref{lemma70}.\thatsall
\begin{lemma}\la{lemma80}There exists a positive constant $C$ such that \bn\ba\la{eq1} \sup_{0\le t\le T}\int_{\O } \te_x^2dx+\int_0^T\int_{\O}  \left( \te_t^2+\te_{xx}^2\right)dxdt\le C .  \ea\en
\end{lemma}

\pf First, multiplying \eqref{2.1} by
$ \te^{\b}\te_t$ and integrating the resultant equality over $\O$, we have
\be\ba\la{8.2} & \int_{\O}  \te^\b\te_t^2dx+ \frac{1}{2} \left(\int_{\O} \frac{(\te^\b\theta_{x})^2}{v}dx\right)_t\\&=- \frac{1}{2}\int_{\O}\frac{(\te^\b\theta_{x})^2u_x}{v^2}dx+\int_{\O} \frac{ \te^{\b }\theta_t\left(-\te u_x+\mu u_x^2+|\nw_x|^2+|\nb_x|^2\right)}{v }dx\\&\le C\sup_{x\in\O} (|u_x|\te^{\b/2})\int_{\O} \te^{3\b/2}\te_x^2dx +\frac{1}{2}\int_{\O} \te^\b\te_t^2dx+C\int_{\O}\te^{\b+2}u_x^2dx\\&\quad+C\int_{\O}\te^\b\left(u_x^4+|\nw_x|^4+|\nb_x|^4\right)dx
\\&\le  C \int_{\O}\te^{2\b}\te_x^2dx\int_{\O}\te^\b\te_x^2dx+\frac{1}{2}\int_{\O}\te^\b\te_t^2dx\\&\quad+C \sup_{x\in\O}(  \te^{2\b+2}+  u_x^4 +|\nw_x|^4+|\nb_x|^4 )+C ,\ea\ee due to
\eqref{6.1}   and \eqref{7.5}.
Combining \eqref{7.1} with H\"{o}lder's inequality gives
 \be\ba\la{8.3} \int_0^T\sup_{x\in\O}u_x^4dt &
 \le C\int_0^T\int_{\O}|u_x^3u_{xx}|dxdt\\&
 \le C\int_0^T\sup_{x\in\O}u_x^2\left(\int_{\O}u_x^2dx\right)^{\frac{1}{2}}\left(\int_{\O}u_{xx}^2dx\right)^{\frac{1}{2}}dt\\&
 \le C\int_0^T\int_{\O}u_{xx}^2dxdt\\&
 \le C.
 \ea\ee
  Using \eqref{6.1} and applying similar arguments to $\nb$ and $\nw$ yields
 \be\ba\la{8.4}\int_0^T\sup_{x\in\O}(|\nb_x|^4+|\nw_x|^4)dt\leq C. \ea\ee
Noticing that for \eqref{2.18}
\be\ba\la{eq9}\sup_{x\in\O} \te^{2\beta+2}&
 \le C\sup_{x\in\O}\left(\int_x^{\infty}\partial_{y}(\te-2)_{+}^{\b+1}dy\right)^2+C\\&
 \le C\int_{(\te>2)(t)}\left(\te^{\b}\te_{x}\right)^2dx+C\\&
 \le C\int_{\O}\frac{\left(\te^{\b}\te_{x}\right)^2}{v}dx+C. \ea\ee
 We  deduce from \eqref{7.5}, \eqref{8.2}--\eqref{eq9}, and   Gr\"onwall's inequality that
\be\ba\la{8.5} \sup_{0 \le t\le T}\int_{\O} \left(\te^\b\theta_{x}\right)^2 dx+\int_0^T\int_{\O} \te^\b\te_t^2dxdt\le C, \ea\ee
which together with \eqref{eq9}  shows
\be\la{8.6}\sup_{(x,t)\in\O\times[0,T]}\te(x,t)\le C.\ee
Thus, both \eqref{2.q27} and \eqref{8.5} lead   to
\be\ba\la{8.7}\sup_{0 \le t\le T}\int_{\O}\theta_{x}^2 dx+\int_0^T\int_{\O} \te_t^2dxdt\le C. \ea\ee

Finally, it follows from \eqref{2.1} that
\bnn\ba \frac{\te^\b\te_{xx}}{v}= -\frac{\b \te^{\b-1}\te_x^2}{v}+\frac{\te^\b\te_x v_x}{v^2}- \frac{\mu u_x^2+|\nb_x|^2+|\nw_x|^2}{v}+ \frac{ \te u_x}{v}+\te_t,\ea\enn
which together with \eqref{2.12},  \eqref{2.q27},   \eqref{4.1}, \eqref{6.1}, \eqref{7.1}, \eqref{8.7}  and   \eqref{8.6} yields
\bnn\ba\la{8.8}\int_0^T\int_{\O} \te_{xx}^2dxdt  & \le C\int_0^T\int_{\O} \left(\te_x^4+\te_x^2v_x^2+u_x^4+|\nb_x|^4+|\nw_x|^4+u_x^2+\te_t^2\right)dxdt\\&
\le C+\int_0^T\left(\sup_{x\in\O}|u_x|^2+\sup_{x\in\O}|\nb_x|^2+\sup_{x\in\O}|\nw_x|^2\right)dt\\&
\quad + C\int_0^T\sup_{x\in\O}\te_x^2 dt\\&
\le C+\frac12\int_0^T \int_{\O} \te_{xx}^2 dx   dt.\ea\enn
 Combining this with \eqref{8.7}  proves \eqref{eq1}
and  finishes the proof of \lemref{lemma80}.\thatsall


\begin{thebibliography}{99}


\bi{az1} Amosov, A. A.,  Zlotnik, A. A.: A difference scheme on a non-uniform mesh for the equations of one-dimensional
 magnetic  gas dynamics. U.S.S.R. Compu. Maths. Math. Phys.,  {\bf 29} (1989), 129--139.

\bi{anto}  Antontsev, S. N., Kazhikhov, A. V., Monakhov, V. N.: Boundary Value Problems
in Mechanics of Nonhomogeneous Fluids. Amsterdam, New York: North-Holland,
1990.

\bi{q1} Cabannes, H.: Theoretical Magnetofluiddynamics. Academic Press, New York (1970)
\bi{cw2} Chen, G. Q., Wang, D. H.: Global solutions for nonlinear magnetohydrodynamics with large initial data. J. Differ. Equ., 182  (2002), 344-376.
\bi{cw1}
Chen, G. Q., Wang, D. H.: Existence and continuous dependence of large solutions for the magnetohydrodynamics
equations. Z. Angew. Math. Phys.,  {\bf 54}  (2003), 608-632.

\bi{q2} Duan, R., Jiang, F., Jiang, S.: On the Rayleigh Taylor instability for incompressible, inviscid magnetohydrodynamic
flows. SIAM J. Appl. Math., {\bf 71} (2011) 1990-2013.

 \bi{fhl1}
Fan, J. S., Huang, S. X., Li, F. C.: Global strong solutions to the planar compressible magnetohydrodynamic equations with large initial data and vaccum.  Kinetic \& Related Models,  {\bf 10} (2017),  1035-1053.

\bi{fjn1}
Fan, J. S., Jiang, S., Nakamura, G.: Vanishing shear viscosity limit in the magnetohydrodynamics equations. Commun. Math. Phys., {\bf 270} (2007), 691-708.


 \bibitem {6} Hoff, D., Tsyganov, E.: Uniqueness and continuous dependence of weak solutions in compressible magnetohydrodynamics.
Z. Angew. Math. Phys., {\bf 56} (2005), 791-804.
\bibitem {7}
 Hu, Y., Ju, Q.: Global large solutions of magnetohydrodynamics with temperature-dependent
heat conductivity. Z. Angew. Math. Phys., $\mathbf{66}$ (2015), 865-889.

\bi{hs1} Huang, B., Shi X. D.:
Nonlinearly exponential stability  of
compressible Navier-Stokes system  with
degenerate heat-conductivity.  J. Differ. Equ., in press. doi.org/10.1016/j.jde.2019.09.006

\bi{hss1}Huang, B., Shi X. D., Sun, Y.: Global strong solutions to   magnetohydrodynamics  with density-dependent viscosity
and degenerate heat-conductivity.     Nonlinearity, $\mathbf{32}$ (2019), 4395-4412.


\bi{q3} Jeffrey, A., Taniuti, T.: Non-Linear Wave Propagation. With Applications to Physics and Magnetohydrodynamics. Academic
Press, New York (1964)

\bibitem {24}
Jenssen, H. K.,  Karper, T. K.: One-dimensional compressible flow with temperature dependent transport coefficients. SIAM Journal on Mathematical Analysis, ${\bf  42}$ (2010),
904-930.
\bi{q4} Jiang, F., Jiang, S., Wang, Y. J.: On the Rayleigh-Taylor instability for incompressible viscous magnetohydrodynamic
equations. Commun. Partial Differ. Equ., {\bf 39} (2014), 399-438.





\bibitem {10}
 Kawashima, S., Nishida, T.: Global solutions to the initial value problem for the equations of onedimensional motion of viscous polytropic gases. J. Math. Kyoto Univ., $\mathbf{21}$ (1981),
825-837.

\bi{ka1} Kazhikhov, A. V.: A  priori  estimates for the solutions of equations of magnetic gas dynamics, Boundary value problems for equations of mathematical physics, Krasnoyarsk, 1987. In Russian

\bibitem {9}
Kazhikhov, A. V., Shelukhin, V. V.: Unique global solution with respect to time of
initial boundary value problems for one-dimensional equations of a viscous gas. J. Appl.
Math. Mech., $\mathbf{41}$ (1977), 273-282.
\bi{q5} Kulikovskiy, A. G., Lyubimov, G. A.: Magnetohydrodynamics. Addison-Wesley, Reading (1965)
\bi{q6} Laudau, L. D., Lifshitz, E. M.: Electrodynamics of Continuous Media. 2nd edn. Pergamon, New York (1984)

\bibitem {L-L}
 Li, J., Liang, Z. L. Some uniform estimates and large-time behavior of solutions
to onedimensional compressible Navier-Stokes system in unbounded domains with
large data. Arch. Rat. Mech. Anal. 220 (2016), 1195-1208.

\bi{py}Li, K. X., Shu X. L., Xu X. J.: Global existence of strong solutions to compressible
Navier-Stokes system with degenerate heat conductivity
in unbounded domains. Math Meth Appl Sci. DOI: 10.1002/mma.5969

\bibitem {13}
 Nash, J.: Le problème de Cauchy pour les équations différentielles d'un fluide général. Bull. Soc. Math. France, $\mathbf{90}$ (1962), 487-497.
\bibitem {28}
Pan, R. H., Zhang, W. Z.: Compressible Navier-Stokes equations with temperature dependent heat conductivities. Commun. Math. Sci., $\mathbf{13}$(2015), 401-425.
\bi{q7} Polovin, R. V., Demutskii, V. P.: Fundamentals of Magnetohydrodynamics. Consultants Bureau, New York (1990)
\bi{kazh} Kazhikhov, A. V. Cauchy problem for viscous gas equations. Siberian Math. J. 23
(1982), 44-49.



\bi{tan} Tani, A.: On the first initial-boundary value problem of compressible viscous fluid motion. Publications
of the Research Institute for Mathematical Sciences, {\bf 13} (1977), 193-253.

\bi{vh}
Vol'pert, A. I., Hudjaev, S. I.: On the Cauchy problem for composite systems of nonlinear differential equations. Math.
USSR-Sb., {\bf 16} (1972), 517-544.

\bibitem {25}
 Wang, D. H.: Large solutions to the initial-boundary value problem for planar magnetohydrodynamics. SIAM J. Appl.
Math., $\mathbf{63}$ (2003), 1424-1441.



\end{thebibliography}
 \end{document}